\documentclass[reqno]{amsart}
\usepackage{float}
\usepackage{amsmath}
\usepackage{graphicx}
\usepackage{latexsym}
\usepackage{amsfonts}
\usepackage{amssymb}

\theoremstyle{plain}
\newtheorem{theorem}{Theorem}

\newtheorem{lemma}[theorem]{Lemma}

\theoremstyle{definition}

\newtheorem{remark}[theorem]{Remark}
\newtheorem*{remark*}{Remark}

\begin{document}
\title{Branching processes in random environment die slowly}
\thanks{Supported in part by the Russian Foundation for Basic Research grant
08-01-00078 and by EPSRC grant number EP/D064988/1}
\author[Vatutin]{Vladimir A. Vatutin}
\address{Steklov Mathematical Institute RAS, Gubkin street 8, 19991, Moscow 
\\
Russia}
\email{vatutin@mi.ras.ru}
\author[Kyprianou]{Andreas E. Kyprianou}
\address{The University of Bath, Claverton Down, Bath, BA2 7AY, UK }
\email{a.kyprianou@bath.ac.uk}
\date{\today }
\maketitle

\begin{abstract}
Let \ $Z_{n,}n=0,1,...,$ be a branching process evolving in the random
environment generated by a sequence of iid generating functions $%
f_{0}(s),f_{1}(s),...,$ and let $S_{0}=0,S_{k}=X_{1}+...+X_{k},k\geq 1,$ be
the associated random walk with $X_{i}=\log f_{i-1}^{\prime }(1),$ $\tau
(m,n)$ be the left-most point of minimum of $\left\{ S_{k},k\geq 0\right\} $
on the interval $[m,n],$ and $T=\min \left\{ k:Z_{k}=0\right\} $. Assuming
that the associated random walk satisfies the Doney condition $P\left(
S_{n}>0\right) \rightarrow \rho \in (0,1),n\rightarrow \infty ,$ we prove
(under the quenched approach) conditional limit theorems, as $n\rightarrow
\infty $, for the distribution of $Z_{nt},$ $Z_{\tau (0,nt)},$ and $Z_{\tau
(nt,n)},$ $t\in (0,1),$ given $T=n$. It is shown that the form of the limit
distributions essentially depends on the location of $\tau (0,n)$ with
respect to the point $nt.$
\end{abstract}

\section{\ Introduction}

Recently a number of papers appeared (see, for instance, \cite{AGKV},\cite%
{GK00}, \cite{V01}-\cite{VD2008}) dealing with branching processes in random
environment in which individuals reproduce independently of each other
according to random offspring distributions which vary from one generation
to the other. The present article complements results established in \cite%
{VD}-\cite{VD2008} \ where critical branching processes in random
environment were investigated under the quenched approach. To give a formal
description of the model under consideration we shall spend some time in
this section introducing notation before proceeding to the main results in
the next section.

\bigskip

Let $\Delta $ be the space of probability measures on $\mathbb{N}%
_{0}:=\{0,1,2,...\}$. Equipped with the metric of total variation $\Delta $
becomes a Polish space. Let $\mathbf{K}$ be a random variable taking values
in $\Delta $. An infinite sequence $\mathbf{\bar{K}}=(\mathbf{K}_{0},\mathbf{%
K}_{1},\ldots )$ of i.i.d. copies of $\mathbf{K}$ is said to form a \emph{%
random environment}. A sequence of $\mathbb{N}_{0}$-valued random variables $%
Z_{0},Z_{1},\ldots $ is called a \emph{branching process in the random
environment} $\mathbf{\bar{K}}$, if $Z_{0}$ is independent of $\mathbf{\bar{K%
}}$ and given $\mathbf{\bar{K}}$ the process $Z=(Z_{0},Z_{1},\ldots )$ is a
Markov chain with 
\begin{equation}
\mathcal{L}\left( Z_{n+1}\;|\;Z_{n}=z(n),\,\mathbf{\bar{K}}=(\mathbf{k}_{0},%
\mathbf{k}_{1},\ldots )\right) \ =\ \mathcal{L}\left( \xi _{n1}+\cdots +\xi
_{nz(n)}\right)  \label{transition}
\end{equation}%
for every $n,\,z(n)\in \mathbb{N}_{0}$ and $\mathbf{k}_{0},\mathbf{k}%
_{1},\ldots \in \Delta $, where $\xi _{n1},\xi _{n2},\ldots $ are i.i.d.
random variables with distribution $\mathbf{k}_{n}$. Setting%
\begin{equation*}
f_{n}(s):=\sum_{j=0}^{\infty }\mathbf{k}_{n}(\{j\})s^{j}
\end{equation*}%
one can rewrite (\ref{transition}) as%
\begin{equation*}
\mathbf{E}\left[ s^{Z_{n+1}}|\;Z_{n}=z(n),\,\mathbf{\bar{K}}=(\mathbf{k}_{0},%
\mathbf{k}_{1},\ldots )\right] =\left( f_{n}\left( s\right) \right)
^{z(n)},\,n\geq 0.
\end{equation*}

Let $\Omega =\left\{ \omega \right\} $ be the space of elementary events with%
\begin{equation*}
\omega =\left\{ \mathbf{k}_{0},\mathbf{k}_{1},\ldots ;z(0),z(1),...\right\} ,
\end{equation*}%
$\mathcal{F}$ be the natural $\sigma $-algebra generated by the subsets of $%
\Omega $ and $\mathbf{P}$ be the corresponding probability measure on the $%
\left( \Omega ,\mathcal{F}\right) $. The triple $\left( \Omega ,\mathcal{F},%
\mathbf{P}\right) $ is be our basic probability space. By $\mathcal{F}%
_{n},n\geq 1,$ we denote the projection of $\mathcal{F}$ on 
\begin{equation*}
\Omega _{n}:=\left\{ \omega _{n}=\left( \mathbf{k}_{0},\mathbf{k}_{1},\ldots
,\mathbf{k}_{n-1};z(0),z(1),...,z(n-1)\right) \right\}
\end{equation*}%
and by $\mathbf{P}_{n}$ the projection of $\mathbf{P}$ on $\mathcal{F}_{n}.$

Let 
\begin{equation*}
f\overset{d}{=}f_{0},\ X_{k}:=\ln f_{k-1}^{\prime}\left( 1\right) ,\,\eta
_{k}:=f_{k-1}^{\prime\prime}\left( 1\right) \left( f_{k-1}^{\prime}\left(
1\right) \right) ^{-2},k\in \mathbb{N}=\{1,2,...\};
\end{equation*}%
\begin{equation*}
X\overset{d}{=}X_{1};\,S_{0}:=0,\;S_{k}:=X_{1}+...+X_{k},\,k\geq 1.
\end{equation*}%
The sequence $\left\{ S_{k},k\geq 0\right\} $ is called the associated
random walk of the corresponding branching process in random environment.

Let 
\begin{equation*}
T:=\min \left\{ k:Z_{k}=0\right\}
\end{equation*}%
and $\tau \left( m,n\right) :=\min \{i\in \lbrack m,n]:\;S_{j}\geq
S_{i},\;j=m,m+1,\ldots ,n\}$ be the left--most point of minimum of the
random walk $\left\{ S_{k},k\geq 0\right\} \;$\ on the discrete time
interval $[m,n].$ In particular we shall write $\tau (n):=\tau \left(
0,n\right) $ is the left--most point of minimum of the random walk on the
discrete time interval $[0,n].$

Properties of branching processes in random environment are specified to a
great extent by the properties of the associated random walk. One of the
most important conditions we impose on the characteristics of our branching
process in this respect is the following Doney condition:

\bigskip

\textbf{Assumption $A1$}. There exists a number $0<\rho <1$ such that 
\begin{equation*}
\mathbf{P}(S_{n}>0)\rightarrow \rho \text{ as }n\rightarrow \infty .
\end{equation*}

\bigskip

As it was shown in \cite{Do}, Condition \textbf{$A1$} is equivalent to the
classical Spitzer condition 
\begin{equation*}
\frac{1}{n}\sum_{k=1}^{n}\mathbf{P}(S_{k}>0)\rightarrow \rho \text{ as }%
n\rightarrow \infty .
\end{equation*}%
Recall that Assumption \textbf{$A1$}~implies 
\begin{equation}
n^{-1}\tau \left( n\right) \overset{d}{\rightarrow }\tau ,\quad n\rightarrow
\infty ,  \label{arcsin}
\end{equation}%
where $\tau $ is a random variable distributed according to the generalized
arcsine law with parameter $\rho $ (\cite{S}, Ch. IV, \S\ 20) and the symbol 
$\overset{d}{\rightarrow }$ means convergence in distribution.

\ Vatutin and Dyakonova, using the quenched approach, have proved in \cite%
{VD3} and \cite{VD05} conditional limit theorems given $\left\{ T>n\right\} $
describing the asymptotic behavior, as $n\rightarrow \infty ,$ of the
distribution of the number of particles at moments $Z_{nt}$ and $Z_{\tau
(nt)},t\in (0,1),~$\ in the branching processses in random environment
respecting Assumption $A1$.

In the present paper we consider the conditioning $\left\{ T=n\right\} $
and, under Assumption A1 and the quenched approach, study the distribution
of the number of particles in our branching process either at moments $nt,$ $%
t\in (0,1),$ or at moments located in a vicinity of points $\tau (nt)$ or $%
\tau (nt,n)$. (On a notational point, here and in the sequel we understand $%
nt$ as $[nt]$, the integer part of $nt$). To formulate our results we need
to specify a number of characteristics related with associated random walks.

Let 
\begin{equation*}
\gamma _{0}:=0,\quad \gamma _{j+1}:=\min (n>\gamma _{j}:S_{n}<S_{\gamma
_{j}})
\end{equation*}%
and 
\begin{equation*}
\Gamma _{0}:=0,\quad \Gamma _{j+1}:=\min (n>\Gamma _{j}:S_{n}>S_{\Gamma
_{j}}),\,j\geq 0,
\end{equation*}%
be the strict descending and strict ascending ladder epochs of $%
\{S_{n},n\geq 0\}$. Introduce the functions 
\begin{eqnarray*}
V(x) &:&=\sum_{j=0}^{\infty }\mathbf{P}(S_{\gamma _{j}}\geq -x),\quad
x>0,\quad V\left( 0\right) =1,\quad V\left( x\right) =0,\quad x<0, \\
U(x) &:&=1+\sum_{j=1}^{\infty }\mathbf{P}(S_{\Gamma _{j}}<x),\quad x>0,\quad
U\left( 0\right) =1,\quad U(x)=0,\quad x<0,
\end{eqnarray*}%
and set%
\begin{equation*}
\Theta (a):=\frac{\sum_{j=a}^{\infty }j^{2}\mathbf{K}{(}\left\{ {j}\right\} {%
)}}{\left( \sum_{r=0}^{\infty }r\mathbf{K}{(}\left\{ {r}\right\} {)}\right)
^{2}},\quad a\in \mathbb{N}.
\end{equation*}

\bigskip

\textbf{Assumption $A2$}. There exist $\varepsilon _{0}>0$ and $a\in \mathbb{%
N}_{0}$ such that 
\begin{eqnarray}
&&\mathbf{E}(\log ^{+}\Theta (a))^{\frac{1}{\rho }+\varepsilon _{0}}<\infty
\quad \mbox{and}\quad \mathbf{E}[V(X)(\log ^{+}\Theta (a))^{1+\varepsilon
_{0}}]<\infty ,  \label{A21} \\
&&\mathbf{E}(\log ^{+}\Theta (a))^{\frac{1}{1-\rho }+\varepsilon
_{0}}<\infty \quad \mbox{and}\quad \mathbf{E}[U(-X)(\log ^{+}\Theta
(a))^{1+\varepsilon _{0}}]<\infty .  \label{A22}
\end{eqnarray}

\bigskip

One can find in \cite{AGKV} and \cite{VD3} more details demonstrating the
importance of these conditions.

Finally, we impose a (rather specific) condition on the form of the
probability generating functions of the underlying branching process in
random environment.

\bigskip

\textbf{Assumption} $A3$. The random offspring generating functions $%
f_{n}(s),\,n=0,1,...,$ are fractional-linear, i.e., they have with
probability 1 the following form 
\begin{equation}
f_{n}(s)=r_{n}+(1-r_{n})\frac{q_{n}}{1-p_{n}s}  \label{linear}
\end{equation}%
where $p_{n}+q_{n}=1,\ p_{n}q_{n}>0.$

\bigskip

It turns out that, many (but not all) of the forthcoming results in this
paper can be proved for the branching processes in random environment
respecting Assumptions $A1-A2$ only. However, the main advantage we gain by
imposing Assumption $A3$ is that if we let%
\begin{eqnarray}
f_{k,n}\left( s\right) &:= &f_{k}(f_{k+1}(...(f_{n-1}\left( s\right)
)...)),\;0\leq k\leq n-1,\;f_{n,n}\left( s\right) :=s,  \notag \\
f_{n,m}\left( s\right) &:= &f_{n-1}(f_{n-2}(...(f_{m}\left( s\right)
)...)),\;n\geq m+1,  \notag
\end{eqnarray}%
and denote 
\begin{equation*}
\eta _{j+1}:=\frac{f_{j}^{\prime \prime }\left( 1\right) }{\left(
f_{j}^{\prime }\left( 1\right) \right) ^{2}},\,b_{m}:=\frac{1}{2}%
\sum_{j=0}^{m-1}\eta _{j+1}e^{-S_{j}}
\end{equation*}%
(noting that both are positive quantities), then Assumption $A3$ implies
(see, for instance, \cite{GK00})that for $m=0,1,....$ and $s\in[0,1]$, 
\begin{equation}
\frac{1}{1-f_{0,m}\left( s\right) }=\frac{e^{-S_{m}}}{1-s}+b_m.
\label{FracExp}
\end{equation}

Later on we need various probability measures specified on the measurable
space $\left( \Omega ,\mathcal{F}\right) .$ To distinguish them we use the
symbols $\mathbf{E}$ and $\mathbf{P}$ to denote the expectation and
probability law generated by the initial measure on the tuples 
\begin{equation*}
(f_{0},f_{1},...,f_{n},\ldots ;Z_{0},Z_{1},...,Z_{n},\ldots )
\end{equation*}%
(there may be some occasional mild abuse of this notation however we do not
anticipate that it will lead to confusion) and the symbols $\mathbf{E}_{\;%
\mathbf{k}},\;\mathbf{P}_{\;\mathbf{k}}\;$ to denote the expectation and
probability law under the fixed environment \ $\mathbf{\bar{k}}=\left( 
\mathbf{k}_{0},\mathbf{k}_{1},...,\mathbf{k}_{n},...\right) $. Along with
the basic probability space $\left( \Omega ,\mathcal{F},\mathbf{P}\right) $
we deal with two its copies $\Big(\Omega ^{-},\mathcal{F}^{-},\mathbf{P}^{-}%
\Big)$ and $\Big(\Omega ^{+},\mathcal{F}^{+},\mathbf{P}^{+}\Big)$. Denote by 
$\left\{ f_{n}^{-},n\geq 0\right\} $ and $\left\{ f_{n}^{+},n\geq 0\right\} $
two sequences of the random environment and by $\left\{ S_{n}^{-},n\geq
0\right\} $ and $\left\{ S_{n}^{+},n\geq 0\right\} $ the corresponding
associate random walks specified on $\Big(\Omega ^{-},\mathcal{F}^{-},%
\mathbf{P}^{-}\Big)$ and $\Big(\Omega ^{+},\mathcal{F}^{+},\mathbf{P}^{+}%
\Big),$ respectively. Later on any characteristics or random variables
related with $\left\{ f_{n}^{-},n\geq 0\right\} $ and $\left\{
f_{n}^{+},n\geq 0\right\} $ are superscripted with the symbols $-$ or $+$,
respectively. Following this practice, we write%
\begin{equation*}
\Gamma ^{-}=\min \{n\geq 1:S_{n}^{-}\geq 0\}
\end{equation*}%
and 
\begin{equation*}
\gamma ^{+}=\min \{n\geq 1:S_{n}^{+}<0\}.
\end{equation*}%
We also study various properties of the pair of branching processes in
random environment given the event 
\begin{equation}
\mathcal{A}_{k,p}:=\left\{ \Gamma ^{-}>k,\gamma ^{+}>p\right\} .
\label{doublegamma}
\end{equation}

Set $D=\sum_{j=1}^{\infty }\mathbf{P}(S_{j}=0)$. In addition to the measures 
$\mathbf{P}^{-}$ and $\mathbf{P}^{+}$ we define measures $\mathbf{\hat{P}}%
^{-} $ and $\mathbf{\hat{P}}^{+}$ on $\Big(\Omega ^{-},\mathcal{F}^{-}\Big)$
and $\Big(\Omega ^{+},\mathcal{F}^{+}\Big)$ $\ $whose restrictions $\mathbf{%
\hat{P}}_{k}^{-}$ and $\mathbf{\hat{P}}_{l}^{+}$ on the $\sigma $-algebras $%
\mathcal{F}_{k}^{-}$ and $\mathcal{F}_{p}^{+},\,k,p\in \mathbb{N}$ are
specified by 
\begin{eqnarray}
\mathbf{\hat{P}}_{k}^{-}(\mathcal{A}^{-}) &=&e^{D}\int_{\mathcal{A}%
^{-}}U\left( -S_{k}^{-}\right) I\{\Gamma ^{-}>k\}d\mathbf{P}^{-},\,\mathcal{A%
}^{-}\in \mathcal{F}_{k}^{-},  \label{defMes-} \\
\mathbf{\hat{P}}_{p}^{+}(\mathcal{A}^{+}) &=&\int_{\mathcal{A}%
^{+}}V(S_{p}^{+})I\{\gamma ^{+}>p\}d\mathbf{P}^{+},\,\mathcal{A}^{+}\in 
\mathcal{F}_{p}^{+}.  \label{defMes+}
\end{eqnarray}%
One can check (see, \cite{VD05}) that the sequences $\left\{ \mathbf{\hat{P}}%
_{k}^{-},\,k\in \mathbb{N}\right\} ,$ $\left\{ \mathbf{\hat{P}}%
_{p}^{-},\,p\in \mathbb{N}\right\} $ consist of well-defined and consistent
probability measures. The probabilistic sense of these measures is rather
transparent: the restriction of $\mathbf{\hat{P}}^{-}$ to $\mathcal{F}%
_{k}^{-}$ is concentrated only on the realizations of the environment whose
associated random walks are negative for the first $k$ steps (except the
starting point) while the restriction of $\mathbf{\hat{P}}^{+}$ to $\mathcal{%
F}_{p}^{+}$ is concentrated only on the realizations of the environment
whose associated random walks are nonnegative for the first $p$ steps.
Indeed, in an appropriate sense, $\mathbf{\hat{P}}^+$ and $\mathbf{\hat{P}}%
^- $ may be thought of as the random walks $S^+$ and $S^-$ conditioned to
stay positive. See for example \cite{BD}.

Further, on the measurable space $\Big(\Omega ^{-}\times \Omega ^{+},%
\mathcal{F}^{-}\times \mathcal{F}^{+}\Big)$ we specify the probability
measure $\mathbf{\hat{P}}:=\mathbf{\hat{P}}^{-}\times \mathbf{\hat{P}}^{+},$
whose projection on the elements of the $\sigma $-algebra $\mathcal{F}%
_{k}^{-}\times \mathcal{F}_{p}^{+}$ is given by 
\begin{equation}
\mathbf{\hat{P}}(\mathcal{A})=e^{D}\int_{\mathcal{A}}U\left(
-S_{k}^{-}\right) V(S_{p}^{+})I\{\mathcal{A}_{k,p}\}d(\mathbf{P}^{-}\times 
\mathbf{P}^{+}),\mathcal{A}\in \mathcal{F}_{k}^{-}\times \mathcal{F}_{p}^{+}
\label{defMes}
\end{equation}%
(see \cite{VD05} for more detailed description of this measure).

\bigskip

With the notation above in hand we list for further references some results
established in \ \cite{VD3} before moving to our main results.

\begin{description}
\item[1)] Under Assumptions $A1-A2$ for any $R\in \mathbb{N}_{0}$ there
exists the limit%
\begin{equation}
q_{R}^{+}:=\lim_{n\rightarrow \infty }f_{R,n}^{+}(0)<1\ \ \hat{\mathbf{P}}%
^{+}\text{- a.s.}  \label{asq}
\end{equation}%
(later on we write for brevity $q^{+}$ for $q_{0}^{+}$);

\item[2)] the tuple of random functions 
\begin{equation}
\zeta _{l,m}^{-}(s):=\frac{1-f_{l,m}^{-}(s)}{e^{S_{l}^{-}-S_{m}^{-}}},\ m\in 
\mathbb{N}_{0},l\geq m+1,  \label{zetam}
\end{equation}%
is such that $\hat{\mathbf{P}}-$ a.s. the limit%
\begin{equation}
\zeta _{\infty ,m}^{-}(s):=\lim_{l\rightarrow \infty }\zeta _{l,m}^{-}(s)
\label{defzetam}
\end{equation}%
exists and is positive and less than 1 for any $s\in \lbrack 0,1).$

For brevity we set $\zeta _{l}^{-}(s):=\zeta _{l,0}^{-}(s)$ and $\zeta
^{-}(s):=\zeta _{\infty ,0}^{-}(s).$ Observe that 
\begin{equation}
\zeta :=\lim_{\min (l,n-l)\rightarrow \infty }\zeta
_{l}^{-}(f_{0,n-l}^{+}(0))=\zeta ^{-}(q^{+})  \label{aszet}
\end{equation}%
exists $\hat{\mathbf{P}}-$ a.s. and, moreover, $\zeta \in (0,1]$ with
probability 1.
\end{description}

\section{Main results}

Now we are ready to formulate the main results of the present paper. Below
our two main theorems we offer some intuition as to their interpretation.

\begin{theorem}
\label{Tclose}Suppose that $A1-$ $A3$ hold. Then for any $R\in \mathbb{Z}$,
any $t\in (0,1)$ and $s\in (0,1]$

1) 
\begin{equation*}
\left\{ \mathbf{E}_{\mathbf{k}}\left[ s^{Z_{\tau (nt)+R}}|T=n\right]
\left\vert \tau (n)\geq nt\right. \right\} \overset{d}{\rightarrow }s\left( 
\frac{1-\Theta _{R}}{1-\Theta _{R}s}\right) ^{2},\,n\rightarrow \infty ,
\end{equation*}%
where%
\begin{equation*}
\Theta _{R}=\left\{ 
\begin{array}{ccc}
\zeta ^{-}\left( f_{0,R}^{+}(0)\right) e^{-S_{R}^{+}} & \text{if} & R\geq 0,
\\ 
&  &  \\ 
\zeta _{\infty ,R}^{-}(0)e^{-S_{R}^{-}} & \text{if} & R<0,%
\end{array}%
\right.
\end{equation*}%
and $\Theta _{R}\in (0,1)$ with probability 1;

2)%
\begin{equation*}
\left\{ \mathbf{E}_{\mathbf{k}}\left[ s^{Z_{\tau (nt,n)+R}}|T=n\right]
\left\vert \tau (n)<nt\right. \right\} \overset{d}{\rightarrow }s\left( 
\frac{1-\theta _{R}}{1-\theta _{R}s}\right) ^{2},\,n\rightarrow \infty ,
\end{equation*}%
where%
\begin{equation*}
\theta _{R}=\left\{ 
\begin{array}{ccc}
q_{R}^{+} & \text{if} & R\geq 0, \\ 
&  &  \\ 
f_{R,0}^{-}(q^{+}) & \text{if} & R<0,%
\end{array}%
\right.
\end{equation*}%
and $\theta _{R}\in (0,1)$ with probability 1.
\end{theorem}

The next theorem deals with the distribution of the number of particles at
moments $nt$, $0<t<1.$ Let 
\begin{equation}
O_{m,n}:=\frac{1-f_{0,n}(0)}{1-f_{m,n}(0)}f_{m,n}(0)b_{m}.  \label{DefOO}
\end{equation}

\begin{theorem}
\label{Tfar}Suppose that Assumptions $A1-$ $A3$ hold. Then, for any $t\in
(0,1)$ and $\lambda \in (0,\infty )$%
\begin{equation*}
\mathbf{E}_{\mathbf{k}}\left[ \exp \left\{ -\lambda \frac{Z_{nt}}{O_{nt,n}}%
\right\} |T=n\right] \overset{p}{\rightarrow }\frac{1}{\left( 1+\lambda
\right) ^{2}}\text{ \ \ as \ }n\rightarrow \infty .
\end{equation*}
\end{theorem}

\bigskip

It is necessary to note that inspite of the unique form of the limit in the
cases $\left\{ \tau (n)\geq nt\right\} $ and $\left\{ \tau (n)<nt\right\} ,$
the behavior of the scaling function $O_{nt,n}$ as $n\rightarrow \infty $ is
different in the two theorems. Later on in Lemmas \ref{LOleftRight} and \ref%
{Comm1} (see also Remark \ref{roughlyspeaking}) we shall see that 
\begin{equation*}
O_{nt,n}I\left\{ \tau (n)\geq nt\right\} \asymp e^{S_{nt}-S_{\tau
(nt)}}I\left\{ \tau (n)\geq nt\right\} ,
\end{equation*}%
i.e., the normalization is, essentially, specified by the past behavior of
the associated random walk. In Theorem \ref{Tfar} 
\begin{equation*}
O_{nt,n}I\left\{ \tau (n)<nt\right\} \asymp e^{S_{nt}-S_{\tau
(nt,n)}}I\left\{ \tau (n)<nt\right\}
\end{equation*}%
i.e., is, essentially, specified by the future behavior of the associated
random walk. 

This fact allows us to give the following non-rigorous interpretation of our
results. If the process dies out at a distant moment $T=n$ then it happens
not as a unique catastrofic event. Before the extinction moment the
evolution of the process consists of a number of "bad" periods where the
size of the population is small. According to Theorem \ref{Tclose}, such
periods are located in the vicinities of \textit{random} points $\tau
(nt)I\left\{ \tau (n)\geq nt\right\} $ and $\tau (nt,n)I\left\{ \tau
(n)<nt\right\} .$ On the other hand, at \textit{nonrandom }points $nt,\,t\in
(0,1),$ the size of the population is, by Theorem \ref{Tfar}, big. For
instance, given $VarX<\infty $, $\log Z_{nt}$ is proportional to $%
S_{nt}-S_{\tau (nt)}\asymp \sqrt{n}$ if $\tau (n)>nt$ and to $S_{nt}-S_{\tau
(nt,n)}\asymp \sqrt{n}$ if $\tau (n)<nt$.. \ Thus, the process dies by
passing through a number of bottlenecks and favorable periods.

\bigskip

The remainder of the paper consists of two sections. The first deals with
the Proof of Theorem \ref{Tclose} and the second with the proof of Theorem %
\ref{Tfar}.

\section{Proof of Theorem \protect\ref{Tclose}}

The proof of Theorem \ref{Tclose} is given right at the very end of this
section. We must first pass through a large number of technical results.

We begin by setting 
\begin{equation*}
\Delta _{m,n}:=\left( \frac{1-f_{m,n-1}(0)}{1-f_{m,n}(0)}\right) ^{2}\left( 
\frac{1-f_{0,n}(0)}{1-f_{0,n-1}(0)}\right) ^{2} = \frac{O^2_{m,n}}{%
O^2_{m,n-1}}\frac{f^2_{m,m-1}(0)}{f^2_{m,n}(0)}
\end{equation*}%
and justifying the following key estimate.

\begin{lemma}
\label{BasicINeq}Under Assumption $A3$ for any $1\leq m<n$%
\begin{equation*}
\frac{1}{1+(1-s)O_{m,n-1}}\Delta _{m,n}^{-1}\leq \mathbf{E}_{\mathbf{k}}%
\left[ s^{Z_{m}}|T=n\right] \leq \frac{1}{1+(1-s)O_{m,n}}\Delta _{m,n}.
\end{equation*}
\end{lemma}

\textbf{Proof}. Clearly,%
\begin{eqnarray*}
\mathbf{E}_{\mathbf{k}}\left[ s^{Z_{m}}|T=n\right] &=&\frac{\mathbf{E}_{%
\mathbf{k}}\left[ s^{Z_{m}};T=n\right] }{\mathbf{P}\left( T=n\right) } \\
&=&\frac{\mathbf{E}_{\mathbf{k}}\left[ s^{Z_{m}};Z_{n}=0\right] -\mathbf{E}_{%
\mathbf{k}}\left[ s^{Z_{m}};Z_{n-1}=0\right] }{f_{0,n}(0)-f_{0,n-1}(0)} \\
&=&\frac{f_{0,m}(sf_{m,n}(0))-f_{0,m}(sf_{m,n-1}(0))}{%
f_{0,m}(f_{m,n}(0))-f_{0,m}(f_{m,n-1}(0))}.
\end{eqnarray*}%
Hence, using the Mean Value Theorem and the monotonicity properties of $%
f_{0,m}^{\prime }(s)$ in $s$ and $f_{m,N}(0)$ in $N$ we get%
\begin{equation*}
s\frac{f_{0,m}^{\prime }(sf_{m,n-1}(0))}{f_{0,m}^{\prime }(f_{m,n}(0))}\leq 
\mathbf{E}_{\mathbf{k}}\left[ s^{Z_{m}}|T=n\right] \leq s\frac{%
f_{0,m}^{\prime }(sf_{m,n}(0))}{f_{0,m}^{\prime }(f_{m,n-1}(0))}.
\end{equation*}%
It is easy to conclude from (\ref{FracExp}) that under Assumption $A3$%
\begin{equation*}
f_{0,m}^{\prime }(s)=\frac{e^{-S_{m}}}{\left( 1-s\right) ^{2}}%
(1-f_{0,m}(s))^{2}.
\end{equation*}%
Therefore,%
\begin{equation*}
\frac{f_{0,m}^{\prime }(sf_{m,n}(0))}{f_{0,m}^{\prime }(f_{m,n-1}(0))}=\frac{%
(1-f_{0,m}(sf_{m,n}(0)))^{2}}{\left( 1-sf_{m,n}(0)\right) ^{2}}\frac{\left(
1-f_{m,n-1}(0)\right) ^{2}}{(1-f_{0,n-1}(0))^{2}}.
\end{equation*}%
Further, again making use of (\ref{FracExp}), we have for $s\in[0,1]$ 
\begin{eqnarray*}
\frac{1-f_{0,m}(sf_{m,n}(0))}{1-sf_{m,n}(0)} &=&\frac{1}{e^{-S_{m}}+\left(
1-sf_{m,n}(0)\right) b_{m}} \\
&=&\frac{1}{e^{-S_{m}}+\left( 1-f_{m,n}(0)\right) b_{m}+(1-s)f_{m,n}(0)b_{m}}
\\
&=&\frac{1}{(1-f_{m,n}(0))\left( 1-f_{0,n}(0)\right)
^{-1}+(1-s)f_{m,n}(0)b_{m}} \\
&=&\frac{1-f_{0,n}(0)}{1-f_{m,n}(0)}\times \frac{1}{1+(1-s)O_{m,n}}
\end{eqnarray*}%
where the third equality follows from the first equality when $s=1$. As a
result we have%
\begin{equation*}
\frac{f_{0,m}^{\prime }(sf_{m,n}(0))}{f_{0,m}^{\prime }(f_{m,n-1}(0))}=\frac{%
1}{1+(1-s)O_{m,n}}\Delta _{m,n}.
\end{equation*}%
Similarly,%
\begin{equation*}
\frac{f_{0,m}^{\prime }(sf_{m,n-1}(0))}{f_{0,m}^{\prime }(f_{m,n}(0))}=\frac{%
1}{1+(1-s)O_{m,n-1}}\Delta _{m,n}^{-1}.
\end{equation*}%
The lemma is proved.\hfill$\square$

\bigskip

To proceed further we need to formulate for future reference several known
statements. In the next two lemmas recall that $\mathcal{A}_{n,r}$ was
defined in (\ref{doublegamma}).

\begin{lemma}
\label{Lfran} (\cite{VD3}, Lemma 3) Let Assumption $A1$ be valid and let $%
T_{l,p},\,l,p\in \mathbb{N},$ be a tuple of uniformly bounded random
variables such that, for any pair $l,p$ the random variable $T_{l,p}$ is
measurable with respect to the $\sigma $-algebra $\mathcal{F}_{l}^{-}\times 
\mathcal{F}_{p}^{+}$ . Then 
\begin{equation}
\lim_{\min (n,r)\rightarrow \infty }\mathbf{E}[T_{l,p}\,|\,\mathcal{A}%
_{n,r}]=\hat{\mathbf{E}}[T_{l,p}].  \label{g11}
\end{equation}%
More generally, if the tuple $\{T_{n,r},\,n,r\in \mathbb{N}\}$ consists of
the uniformly bounded random variables which are adopted to the flow of the $%
\sigma -$algebras\break $\{\mathcal{F}_{n}^{-}\times \mathcal{F}%
_{r}^{+}\}_{n\geq 1,r\geq 1},$ and $\lim_{\min (n,r)\rightarrow \infty
}T_{n,r}=:T$ exists $\hat{\mathbf{P}}-$ a.s. then 
\begin{equation}
\lim_{\min (n,r)\rightarrow \infty }\mathbf{E}[T_{n,r}\,|\,\mathcal{A}%
_{n,r}]=\hat{\mathbf{E}}[T].  \label{g22}
\end{equation}
\end{lemma}

\begin{lemma}
\label{Cfran} (\cite{VD3}, Lemma 4) Let Assumption $A1$ be valid and let $T$
and $T_{l,p},\,l,p\in \mathbb{N}$ be a tuple of random variables meeting the
conditions of Lemma \ref{Lfran}. If $T_{l,p}^{\ast },\,l,p\in \mathbb{N}$ is
a tuple of uniformly bounded random variables such that for any pair $l,p$
the random variable, $T_{l,p}^{\ast }$ is measurable with respect to the $%
\sigma $-algebra $\mathcal{F}_{l}^{-}\times \mathcal{F}_{p}^{+}$ and 
\begin{equation*}
\mathbf{E}[T_{\tau (n),n}^{\ast }\,|\,\tau (n)=l]=\mathbf{E}[T_{l,n-l}\,|\,%
\mathcal{A}_{l,n-l}],
\end{equation*}%
then 
\begin{equation}
\lim_{n\rightarrow \infty }\mathbf{E}[T_{\tau (n),n}^{\ast }]=\hat{\mathbf{E}%
}[T].  \label{g112}
\end{equation}
\end{lemma}

\textbf{\ }

For $0\leq m\leq n$ set 
\begin{equation}
\alpha _{n}(m):=\frac{1-f_{0,m}(0)}{1-f_{0,n}(0)},\qquad \beta _{n}(m):=%
\frac{1-f_{0,n}\left( 0\right) }{e^{S_{m}}\left( 1-f_{m,n}\left( 0\right)
\right) }.\,  \label{alfbet}
\end{equation}%
Clearly, $\alpha _{n}(m)\geq 1,\ \beta _{n}(m)\leq 1,$ and, in addition,%
\begin{eqnarray}
O_{m,n} &=&\frac{1}{\alpha _{n}(m)}\times \frac{f_{m,n}(0)}{1-f_{m,n}(0)}%
-f_{m,n}(0)\beta _{n}(m)  \notag \\
&=&f_{m,n}(0)\beta _{n}(m)\left( \frac{e^{S_{m}}}{1-f_{0,m}(0)}-1\right) .
\label{DDoubleRepres}
\end{eqnarray}

\begin{lemma}
\label{min1} (\cite{VD3}, Lemma 16) Assume that $A1-$ $A2$ hold. Then for $%
R\in \mathbb{Z}$ and any $\varepsilon >0$ 
\begin{equation*}
\limsup_{R,n\rightarrow \infty }\mathbf{P}\left( \alpha _{n}(\tau
(n)+R)>1+\varepsilon \right) =0.
\end{equation*}
\end{lemma}

\begin{remark}
Since the random variable $\alpha _{n}(\tau (n)+R)$ is not defined for $\tau
(n)+R<0$ or $\tau (n)+R>n$ we should formally write the statement of the
lemma as 
\begin{equation*}
\limsup_{R,n\rightarrow \infty }\mathbf{P}\left( \alpha _{n}(\tau
(n)+R)>1+\varepsilon ;n\geq \tau (n)+R\geq 0\right) =0
\end{equation*}%
However, by the generalized arcsine law for each fixed $R\in Z$%
\begin{equation}
\lim_{n\rightarrow \infty }\mathbf{P}\left( \tau (n)+R\notin \lbrack
0,n]\right) =0.  \label{Inter1}
\end{equation}%
For this reason here and in what follows we agree to treat $\alpha _{n}(\tau
(n)+R)$ as $\alpha _{n}(0)$ if $\tau (n)+R<0$ and $\alpha _{n}(n)$ if $\tau
(n)+R>n.$ Similar agreement will be kept for other functions which involve $%
\tau (n)+R,$\textbf{\ }$\tau (nt)+R$\textbf{\ or }$\tau (nt,n)+R.$
\end{remark}

\begin{lemma}
\label{min2} (\cite{VD3}, Corollary 3) Assume that $A1-$ $A2$ hold. For any $%
t\in \left( 0,1\right] $ and $\varepsilon >0$ 
\begin{equation*}
\limsup_{n\rightarrow \infty }\mathbf{P}\left( \,\alpha
_{n}(nt)>1+\varepsilon \,\left\vert \,\tau \left( n\right) <nt\right.
\right) =0.
\end{equation*}
\end{lemma}

The next statement complements Lemmas \ref{min1} and \ref{min2}.

\begin{lemma}
\label{Lprospect1} Assume that $A1-$ $A2$ hold. Then for $R\in \mathbb{Z},$
any $t\in \left( 0,1\right] ,$ and any $\varepsilon >0$%
\begin{equation*}
\limsup_{R,n\rightarrow \infty }\mathbf{P}\left( \alpha _{n}(\tau
(nt,n)+R)>1+\varepsilon \left\vert \,\tau \left( n\right) <nt\right. \right)
=0.
\end{equation*}
\end{lemma}

\textbf{Proof}. First we note that \ $\alpha _{n}(m)\downarrow 1$ as $%
m\uparrow n$. Hence 
\begin{eqnarray*}
\alpha _{n}(\tau (nt,n)+R)I\left\{ nt\leq \tau (nt,n)+R\leq n\right\}  &\leq
&\alpha _{n}(nt)I\left\{ nt\leq \tau (nt,n)+R\leq n\right\}  \\
&\leq &\alpha _{n}(nt).
\end{eqnarray*}%
Therefore, 
\begin{eqnarray}
\mathbf{P}\left( \alpha _{n}(\tau (nt,n)+R)>1+\varepsilon ;\tau \left(
n\right) <nt\right)  &\leq &\mathbf{P}\left( \alpha _{n}(nt)>1+\varepsilon
;\tau \left( n\right) <nt\right)   \notag \\
&&+\mathbf{P}\left( \tau (nt,n)+R\notin \lbrack nt,n]\right) .  \label{Adr1}
\end{eqnarray}%
Observing that $\left\{ \tau (nt,n)+R>n\right\} =\oslash $ for $R<0$ and $%
\left\{ \tau (nt,n)+R<nt\right\} =\oslash $ for $R>0$ and recalling the
generalized arcsine law (which, under Assumption $A1,$ holds for $n^{-1}\tau
(nt,n)$ (see \cite{S})) we see that for each fixed $R$ 
\begin{equation}
\lim_{n\rightarrow \infty }\mathbf{P}\left( \tau (nt,n)+R\notin \lbrack
nt,n]\right) =0.  \label{Interm}
\end{equation}%
Passing now to the limit as $n\rightarrow \infty $ in the both sides of (\ref%
{Adr1}) and recalling Lemma \ref{min2} we see that%
\begin{equation*}
\lim \sup_{n\rightarrow \infty }\mathbf{P}\left( \alpha _{n}(\tau
(nt,n)+R)>1+\varepsilon ;\tau \left( n\right) <nt\right) =0.
\end{equation*}%
Hence the statement of the lemma follows. \hfill $\square $



\begin{lemma}
\label{Lprospect2}Assume that $A1-$ $A2$ hold. Then for any $R\in \mathbb{Z}$%
, any $t\in (0,1),$ and $N\in \{n-1,$ $n\}$%
\begin{equation*}
\left\{ f_{\tau (nt,n)+R,N}(0)|\tau (n)<nt\right\} \overset{d}{\rightarrow }%
\theta _{R}\text{ as }n\rightarrow \infty ,
\end{equation*}%
where $\theta _{R}$ is the same as in Theorem \ref{Tclose}.
\end{lemma}

\textbf{Proof}. Consider first $R<0.$ 
Clearly, 
\begin{eqnarray*}
&&\mathbf{E}\left[ e^{-\lambda f_{\tau (nt,n)+R,N}(0)};\tau (n)<nt\right]  \\
&=&\mathbf{E}\left[ e^{-\lambda f_{\tau (nt,n)+R,N}(0)}\right] -\mathbf{E}%
\left[ e^{-\lambda f_{\tau (nt,n)+R,N}(0)};\tau (n)\geq nt\right] .
\end{eqnarray*}%
Introduce two independent environmental sequences $\left\{ f_{n}^{-},n\geq
0\right\} $ and $\left\{ f_{n}^{+},n\geq 0\right\} $ and the respective
associated random walks $\left\{ S_{n}^{-},n\geq 0\right\} $ and $\left\{
S_{n}^{+},n\geq 0\right\} $. In the notation of Lemmas \ref{Lfran} and \ref%
{Cfran} we set 
\begin{equation*}
T_{\tau (n(1-t)),n(1-t)}^{\ast }:=e^{-\lambda f_{\tau
(n(1-t))+R,N-n(1-t)}(0)},\,T_{l,n-l}:=e^{-\lambda
f_{R,0}^{-}(f_{0,N-n(1-t)-l}^{+}(0))}.
\end{equation*}%
Observing that\textbf{\ }%
\begin{equation*}
e^{-\lambda f_{\tau (n(1-t))+R,N-n(1-t)}(0)}\overset{d}{=}e^{-\lambda
f_{\tau (nt,n)+R,N}(0)}
\end{equation*}%
and%
\begin{equation*}
f_{R,0}^{-}(f_{0,N-n(1-t)-l}(0))\rightarrow f_{R,0}^{-}(q^{+})
\end{equation*}%
$\mathbf{\hat{P}}$ $-$a.s. as $\min (l,N-n(1-t)-l)\rightarrow \infty $ we
conclude by Lemma \ref{Cfran} and the generalized arcsine law that 
\begin{eqnarray*}
\mathbf{E}\left[ e^{-\lambda f_{\tau (nt,n)+R,N}(0)}\right]  &=&\mathbf{E}%
\left[ e^{-\lambda f_{\tau (n(1-t))+R,N-n(1-t)}(0)}\right]  \\
&\rightarrow &\mathbf{\hat{E}}\left[ e^{-\lambda f_{R,0}^{-}(q^{+})}\right]
,\,n\rightarrow \infty .
\end{eqnarray*}%
Further we have 
\begin{eqnarray*}
\lefteqn{\mathbf{E}\left[ e^{-\lambda f_{\tau (nt,n)+R,N}(0)};\tau (n)\geq nt%
\right] =\mathbf{E}\left[ e^{-\lambda f_{\tau (nt,n)+R,N}(0)};\tau (n)\geq nt%
\right] } \\
&=&\sum_{j=nt}^{n}\mathbf{E}\left[ e^{-\lambda f_{\tau (n)+R,N}(0)};\tau
(n)=j\right]  \\
&=&\sum_{j=nt}^{n}\mathbf{E}\left[ e^{-\lambda f_{j+R,N}(0)};\min_{0\leq
r<j}S_{r}>S_{j},\min_{j+1\leq r\leq n}S_{r}\geq S_{j}\right]  \\
&=&\sum_{j=nt}^{n}\mathbf{E}\left[ e^{-\lambda
f_{R,0}^{-}(f_{0,N-j}^{+}(0))};A_{j,n-j}\right]  \\
&=&\sum_{j=nt}^{n}\mathbf{E}\left[ e^{-\lambda
f_{R,0}^{-}(f_{0,N-j}^{+}(0))}|A_{j,n-j}\right] \mathbf{P}(\tau (n)=j)
\end{eqnarray*}%
and since $f_{R,0}^{-}(f_{0,N-j}^{+}(0))\rightarrow f_{R,0}^{-}(q^{+})$ $%
\mathbf{\hat{P}}-$a.s as $n-j\rightarrow \infty $ we have by Lemma~\ref%
{Lfran} that 
\begin{equation*}
\lim_{n-j\rightarrow \infty }\mathbf{E}\left[ e^{-\lambda
f_{R,N-j}(0)}|A_{j,n-j}\right] =\mathbf{\hat{E}}\left[ e^{-\lambda
f_{R,0}^{-}(q^{+})}\right] .
\end{equation*}%
This and the generalized arcsine law give%
\begin{equation*}
\lim_{n\rightarrow \infty }\mathbf{E}\left[ e^{-\lambda f_{\tau
(nt,n)+R,N}(0)};\tau (n)\geq nt\right] =\mathbf{\hat{E}}\left[ e^{-\lambda
f_{R,0}^{-}(q^{+})}\right] \mathbf{P}\left( \tau \geq t\right) .
\end{equation*}



Thus,%
\begin{equation*}
\lim_{n\rightarrow \infty }\mathbf{E}\left[ e^{-\lambda f_{\tau
(nt,n)+R,N}(0)I\left\{ \mathbf{\tau }(nt,n)+R\geq 0\right\} };\tau (n)<nt%
\right] =\mathbf{\hat{E}}\left[ e^{-\lambda f_{R,0}^{-}(q^{+})}\right] 
\mathbf{P}\left( \tau \leq t\right)
\end{equation*}%
proving the lemma for $R<0$.

The case $R\geq 0$ can be treated in a similar way by observing that, in
this case, for all $j\leq N-R$ 
\begin{equation*}
\mathbf{E}\left[ e^{-\lambda f_{\tau (n)+R,N}(0)};\tau (n)=j\right] =\mathbf{%
E}\left[ e^{-\lambda f_{R,N-j}^{+}(0))}|A_{j,n-j}\right] \mathbf{P}(\tau
(n)=j)
\end{equation*}%
and that (\ref{asq}) holds.\hfill $\square $

\bigskip

In what follows we need some properties of the random variable $\beta
_{n}(m) $.

\begin{lemma}
\label{min22}(\cite{VD3}, Lemma 17) Assume that $A1-$ $A2$ hold. For any $%
\varepsilon >0$ and $N\in \left\{ n-1,n\right\} $%
\begin{equation}
\limsup_{R,n\rightarrow \infty }\mathbf{P}\left( \beta _{N}\left( \tau
\left( n\right) +R\right) >\varepsilon \right) =0.  \label{betvan}
\end{equation}
\end{lemma}

\begin{lemma}
\label{min4}(\cite{VD3}, Corollary 4) Assume that $A1-$ $A2$ hold. For any $%
t\in \left( 0,1\right) $ and $\varepsilon >0$%
\begin{eqnarray}
&&\limsup_{n\rightarrow \infty }\mathbf{P}\left( \beta _{n}\left( nt\right)
<1-\varepsilon ;\tau \left( n\right) \geq nt\right) =0,  \label{min41} \\
&&\limsup_{n\rightarrow \infty }\mathbf{P}\left( \beta _{n}\left( nt\right)
>\varepsilon ;\tau \left( n\right) <nt\right) =0.  \label{min42}
\end{eqnarray}
\end{lemma}

\begin{lemma}
\label{EstO}Assume that $A1-$ $A2$ hold. For any $R\in \mathbb{Z}$, any $%
t\in (0,1),$ and $N\in \left\{ n-1,n\right\} $ 
\begin{equation*}
\left\{ O_{\tau (nt,n)+R,N}|\tau (n)<nt\right\} \overset{d}{\rightarrow }%
\frac{\theta _{R}}{1-\theta _{R}}.
\end{equation*}
\end{lemma}

\textbf{Proof}. We have%
\begin{equation}
O_{\tau (nt,n)+R,N}=\frac{1}{\alpha _{N}(\tau (nt,n)+R)}\times \frac{f_{\tau
(nt,n)+R,N}(0)}{1-f_{\tau (nt,n)+R,N}(0)}-f_{\tau (nt,n)+R,N}(0)\beta
_{N}(\tau (nt,n)+R).  \label{uu0}
\end{equation}%
By Lemmas \ref{Lprospect1} and \ref{Lprospect2}%
\begin{equation}
\frac{1}{\alpha _{N}(\tau (nt,n)+R)}\times \frac{f_{\tau (nt,n)+R,N}(0)}{%
1-f_{\tau (nt,n)+R,N}(0)}\overset{d}{\rightarrow }\frac{\theta _{R}}{%
1-\theta _{R}},\,n\rightarrow \infty .  \label{uu1}
\end{equation}%
Further, by (\ref{min42})%
\begin{eqnarray}
\lefteqn{\beta _{N}(\tau (nt,n)+R)I\left\{ \tau (n)<nt\leq \tau
(nt,n)+R\right\} } \\
&\leq &\beta _{N}(nt)I\left\{ \tau (n)<nt\leq \tau (nt,n)+R\right\}   \notag
\\
&\leq &\beta _{N}(nt)I\left\{ \tau (n)<nt\right\} \overset{p}{\rightarrow }0
\label{uu2}
\end{eqnarray}%
as $n\rightarrow \infty $. Using (\ref{uu1}) and (\ref{uu2}) to evaluate (%
\ref{uu0}) proves the lemma.\hfill $\square $

\begin{lemma}
\label{Lpreelimit1}Assume that $A1-$ $A2$ hold. For any $t\in \left(
0,1\right) ,$ any fixed $R\in \mathbb{Z}$, and $\varepsilon >0$%
\begin{equation*}
\limsup_{n\rightarrow \infty }\mathbf{P}\left( \beta _{n}\left( \tau
(nt)+R\right) <1-\varepsilon ;\tau \left( n\right) \geq nt\right) =0.
\end{equation*}
\end{lemma}

\textbf{Proof}. Since 
\begin{equation*}
e^{S_{m}}\left( 1-f_{m,n}\left( s\right) \right) \leq e^{S_{m+1}}\left(
1-f_{m+1,n}\left( s\right) \right)
\end{equation*}%
for any $m<n,$ the elements of the sequence $\left\{ \beta _{n}(m),0\leq
m\leq n\right\} $ are monotone decreasing in $m$ for any fixed $n.$ On the
other hand, for any $\delta >0$ there exists $\varepsilon _{1}>0$ such that%
\begin{equation*}
\mathbf{P}\left( \tau (nt)+R>nt(1-\varepsilon _{1})\right) <\delta
\end{equation*}%
for all $n\geq n_{0}=n_{0}(\delta ,\varepsilon _{1}).$ Thus, we have for all 
$n\geq n_{0}$ 
\begin{eqnarray*}
\lefteqn{\mathbf{P}\left( \beta _{n}\left( \tau (nt)+R\right) <1-\varepsilon
;\tau \left( n\right) \geq nt\right) } \\
&\leq &\delta +\mathbf{P}\left( \beta _{n}\left( nt(1-\varepsilon
_{1})\right) <1-\varepsilon ;\tau \left( n\right) \geq nt\right) \\
&\leq &\delta +\mathbf{P}\left( \beta _{n}\left( nt(1-\varepsilon
_{1})\right) >1-\varepsilon ;\tau \left( n\right) \geq nt(1-\varepsilon
_{1})\right) .
\end{eqnarray*}%
To complete the proof it remains to recall (\ref{min41}) .\hfill $\square $

\begin{lemma}
\label{Lprelim3}Assume that $A1-$ $A2$ hold. For any $R\in \mathbb{Z}$, any $%
t\in (0,1),$ and $N\in \left\{ n-1,n\right\} $ 
\begin{equation*}
\left\{ f_{\tau (nt)+R,N}(0)|\tau (n)\geq nt\right\} \overset{d}{\rightarrow 
}1,\,n\rightarrow \infty .
\end{equation*}
\end{lemma}

\textbf{Proof}. For a fixed $\varepsilon \in (0,1)$ introduce the event 
\begin{equation*}
\mathcal{H}_{n}(\varepsilon )=\left\{ \omega :\left\{ \tau (nt),\tau
(nt)+R,\tau (n)\right\} \cap \lbrack nt(1-\varepsilon ),nt(1+\varepsilon
)]=\oslash \right\} 
\end{equation*}%
and let $\mathcal{\bar{H}}_{n}(\varepsilon )$ be the complement of \ $%
\mathcal{H}_{n}(\varepsilon ).$ By the generalized arcsine law 
\begin{eqnarray*}
\lim_{\varepsilon \downarrow 0}\lim_{n\rightarrow \infty }\mathbf{P}\left( 
\mathcal{\bar{H}}_{n}(\varepsilon )\right)  &\leq &\lim_{\varepsilon
\downarrow 0}\lim_{n\rightarrow \infty }\mathbf{P}\left( \tau (nt)\in
\lbrack nt(1-\varepsilon ),nt(1+\varepsilon )]\right)  \\
&&+\lim_{\varepsilon \downarrow 0}\lim_{n\rightarrow \infty }\mathbf{P}%
\left( \tau (nt)+R\in \lbrack nt(1-\varepsilon ),nt(1+\varepsilon )]\right) 
\\
&&+\lim_{\varepsilon \downarrow 0}\lim_{n\rightarrow \infty }\mathbf{P}%
\left( \tau (n)\in \lbrack nt(1-\varepsilon ),nt(1+\varepsilon )]\right)  \\
&=&0.
\end{eqnarray*}%
Hence, to prove the lemma it is suffices to show that 
\begin{equation*}
\left\{ f_{\tau (nt)+R,N}(0)|\mathcal{H}_{n}(\varepsilon );\tau (n)\geq
nt\right\} \overset{d}{\rightarrow }1,\,n\rightarrow \infty .
\end{equation*}%
Clearly,%
\begin{eqnarray*}
\lefteqn{\left( 1-f_{\tau (nt)+R,N}(0)\right) I\left\{ \mathcal{H}%
_{n}(\varepsilon );\tau (n)\geq nt\right\} } \\
&\leq &\left( 1-f_{\tau (nt)+R,\tau (n)}(0)\right) I\left\{ \mathcal{H}%
_{n}(\varepsilon );\tau (n)>nt(1+\varepsilon )\right\}  \\
&\leq &e^{S_{\tau (n)}-S_{\tau (nt)+R}}I\left\{ \mathcal{H}_{n}(\varepsilon
);\tau (n)>nt(1+\varepsilon )\right\} \leq e^{S_{\tau (n)}-S_{\tau
(nt)}}I\left\{ \mathcal{H}_{n}(\varepsilon );\tau (n)>nt(1+\varepsilon
)\right\}  \\
&\leq &e^{S_{\tau (n)}-S_{\tau (nt)}}I\left\{ \tau (n)>nt(1+\varepsilon
)\right\} .
\end{eqnarray*}%
Thus, for any $\varepsilon _{1}\in (0,1)$ 
\begin{eqnarray}
\lefteqn{\mathbf{P}\left( 1-f_{\tau (nt)+R,N}(0)\geq \varepsilon _{1};%
\mathcal{H}_{n}(\varepsilon );\tau (n)\geq nt\right) }  \notag \\
&\leq &\varepsilon _{1}^{-1}\mathbf{E}\left[ 1-f_{\tau (nt)+R,N}(0);\mathcal{%
H}_{n}(\varepsilon );\tau (n)\geq nt\right]   \notag \\
&\leq &\varepsilon _{1}^{-1}\mathbf{E}\left[ e^{S_{\tau (n)}-S_{\tau
(nt)}};\tau (n)>nt(1+\varepsilon )\right] .  \label{ProPrel}
\end{eqnarray}%
Hence, using the notation introduced in the proof of Lemma \ref{Lprospect2}
and the duality principle for random walks it is not difficult to check that%
\begin{eqnarray*}
\lefteqn{\mathbf{E}\left[ e^{S_{\tau (n)}-S_{\tau (nt)}};\tau
(n)>nt(1+\varepsilon )\right] } \\
&=&\sum_{nt(1+\varepsilon )<k\leq n}\mathbf{E}\left[ e^{S_{k}-\min_{0\leq
i\leq nt}S_{i}};\tau (n)=k\right]  \\
&=&\sum_{nt(1+\varepsilon )<k\leq n}\mathbf{E}\left[ e^{\min_{k-nt\leq l\leq
k}S_{l}^{-}};A_{k,n-k}\right]  \\
&\leq &\sum_{nt(1+\varepsilon )<k\leq n}\mathbf{E}\left[ e^{\min_{nt%
\varepsilon \leq l\leq k}S_{l}^{-}}|A_{k,n-k}\right] \mathbf{P}\left( \tau
(n)=k\right) .
\end{eqnarray*}%
Since $S_{l}^{-}\rightarrow -\infty $ $\mathbf{\hat{P}}$-a.s \ as $%
l\rightarrow \infty ,$ we have by Lemma \ref{Lfran} 
\begin{equation*}
\lim_{n\rightarrow \infty }\mathbf{E}\left[ e^{\min_{nt\varepsilon \leq
l\leq k}S_{l}^{-}}|A_{k,n-k}\right] =0.
\end{equation*}%
Now the dominanted convergence theorem gives%
\begin{equation*}
\lim_{n\rightarrow \infty }\mathbf{E}\left[ e^{S_{\tau (n)}-S_{\tau
(nt)}};\tau (n)>nt(1+\varepsilon )\right] =0.
\end{equation*}%
Combining this fact with (\ref{ProPrel}) complete the proof of the
lemma.\hfill $\square $

\bigskip

\begin{lemma}
\label{Lprelim4}Assume that $A1-$ $A2$ hold. For any $R\in \mathbb{Z}$ and
any $t\in (0,1)$%
\begin{equation*}
\left\{ \frac{1-f_{0,\tau (nt)+R}(0)}{e^{S_{\tau (nt)+R}}}|\tau (n)\geq
nt\right\} \overset{d}{\rightarrow }\Theta _{R}
\end{equation*}%
where $\Theta _{R}$ is the same as in Theorem \ref{Tclose}.
\end{lemma}

\textbf{Proof}. Let $R\geq 0$ be fixed. We have%
\begin{eqnarray*}
\lefteqn{\mathbf{E}\left[ \exp \left\{ -\lambda \frac{1-f_{0,\tau (nt)+R}(0)%
}{e^{S_{\tau (nt)+R}}}\right\} ;\tau (n)\geq nt\right] } \\
&=&\mathbf{E}\left[ \exp \left\{ -\lambda \frac{1-f_{0,\tau (nt)+R}(0)}{%
e^{S_{\tau (nt)+R}}}\right\} \right] \\
&&\hspace{2cm}-\mathbf{E}\left[ \exp \left\{ -\lambda \frac{1-f_{0,\tau
(nt)+R}(0)}{e^{S_{\tau (nt)+R}}}\right\} ;\tau (n)<nt\right] .
\end{eqnarray*}%
Introducing once again two independent environmental sequences $\left\{
f_{n}^{-},n\geq 0\right\} $ and $\left\{ f_{n}^{+},n\geq 0\right\} ,$ setting%
\begin{equation*}
T_{\tau (nt),nt}^{\ast }:=\exp \left\{ -\lambda \frac{1-f_{0,\tau (nt)+R}(0)%
}{e^{S_{\tau (nt)+R}}}\right\} ,T_{l,nt-l}:=\exp \left\{ -\lambda \frac{%
1-f_{l,0}^{-}(f_{0,R}^{+}(0))}{e^{S_{l}^{-}}e^{S_{R}^{+}}}\right\}
\end{equation*}%
and observing that%
\begin{equation*}
\frac{1-f_{l,0}^{-}(f_{0,R}^{+}(0))}{e^{S_{l}^{-}}e^{S_{R}^{+}}}\rightarrow 
\frac{\zeta ^{-}(f_{0,R}^{+}(0))}{e^{S_{R}^{+}}}
\end{equation*}%
$\mathbf{\hat{P}}$ a.s. as $\min (l,n(1-t)-l)\rightarrow \infty ,$ we
conclude by Lemma \ref{Cfran} that for any $\lambda \in \left( 0,\infty
\right) $%
\begin{equation*}
\mathbf{E}\left[ \exp \left\{ -\lambda \frac{1-f_{0,\tau (nt)+R}(0)}{%
e^{S_{\tau (nt)+R}}}\right\} \right] \rightarrow \mathbf{\hat{E}}\left[ \exp
\left\{ -\lambda \frac{\zeta ^{-}(f_{0,R}^{+}(0))}{e^{S_{R}^{+}}}\right\} %
\right]
\end{equation*}%
as $n\rightarrow \infty .$ Further, the same as in Lemma \ref{Lprelim3} 
\begin{eqnarray*}
\lefteqn{\mathbf{E}\left[ \exp \left\{ -\lambda \frac{1-f_{0,\tau (nt)+R}(0)%
}{e^{S_{\tau (nt)+R}}}\right\} ;\tau (n)<nt\right] } \\
&=&\mathbf{E}\left[ \exp \left\{ -\lambda \frac{1-f_{0,\tau (n)+R}(0)}{%
e^{S_{\tau (n)+R}}}\right\} ;\tau (n)<nt\right] \\
&=&\sum_{j=0}^{nt-1}\mathbf{E}\left[ \exp \left\{ -\lambda \frac{%
1-f_{0,j+R}(0)}{e^{S_{j+R}}}\right\} ;\tau (n)=j\right] \\
&=&\sum_{j=0}^{nt-1}\mathbf{E}\left[ \exp \left\{ -\lambda \frac{%
1-f_{0,j+R}(0)}{e^{S_{j+R}}}\right\} ;\min_{0\leq
r<j}S_{r}>S_{j},\min_{j+1\leq r\leq n}S_{r}\geq S_{j}\right] \\
&=&\sum_{j=0}^{nt-1}\mathbf{E}\left[ \exp \left\{ -\lambda \frac{%
1-f_{j,0}^{-}(f_{0,R}^{+}(0))}{e^{S_{j}^{-}}e^{S_{R}^{+}}}\right\} ;A_{j,n-j}%
\right] \\
&=&\sum_{j=0}^{nt-1}\mathbf{E}\left[ \exp \left\{ -\lambda \frac{%
1-f_{j,0}^{-}(f_{0,R}^{+}(0))}{e^{S_{j}^{-}}e^{S_{R}^{+}}}\right\}
\left\vert A_{j,n-j}\right. \right] \mathbf{P}\left( \tau (n)=j\right) .
\end{eqnarray*}%
By\ (\ref{aszet}) 
\begin{equation*}
\frac{1-f_{j,0}^{-}(f_{0,R}^{+}(0))}{e^{S_{j}^{-}}e^{S_{R}^{+}}}\rightarrow 
\frac{\zeta ^{-}(f_{0,R}^{+}(0))}{e^{S_{R}^{+}}}
\end{equation*}%
$\mathbf{\hat{P}}$ a.s. as $\min (j,n-j)\rightarrow \infty .$ Hence we see
by Lemma \ref{Lfran} that 
\begin{equation*}
\lim_{\min (j,n-j)\rightarrow \infty }\mathbf{E}\left[ \exp \left\{ -\lambda 
\frac{\zeta ^{-}(f_{0,R}^{+}(0))}{e^{S_{R}^{+}}}\right\} \left\vert
A_{j,n-j}\right. \right] =\mathbf{\hat{E}}\left[ \exp \left\{ -\lambda \frac{%
\zeta ^{-}(f_{0,R}^{+}(0))}{e^{S_{R}^{+}}}\right\} \right] .
\end{equation*}%
This and the generalized arcsine law give%
\begin{eqnarray*}
\lefteqn{\lim_{n\rightarrow \infty }\mathbf{E}\left[ \exp \left\{ -\lambda 
\frac{1-f_{0,\tau (nt)+R}(0)}{e^{S_{\tau (nt)+R}}}\right\} ;\tau (n)<nt%
\right] } \\
&=&\mathbf{\hat{E}}\left[ \exp \left\{ -\lambda \frac{\zeta
^{-}(f_{0,R}^{+}(0))}{e^{S_{R}^{+}}}\right\} \right] \mathbf{P}\left( \tau
\leq t\right) .
\end{eqnarray*}%
Thus,%
\begin{eqnarray*}
\lefteqn{\lim_{n\rightarrow \infty }\mathbf{E}\left[ \exp \left\{ -\lambda 
\frac{1-f_{0,\tau (nt)+R}(0)}{e^{S_{\tau (nt)+R}}}\right\} ;\tau (n)\geq nt%
\right] } \\
&=&\mathbf{\hat{E}}\left[ \exp \left\{ -\lambda \frac{\zeta
^{-}(f_{0,R}^{+}(0))}{e^{S_{R}^{+}}}\right\} \right] \mathbf{P}\left( \tau
\geq t\right)
\end{eqnarray*}%
and the statement of the lemma for $R\geq 0$ follows. The case $R<0$ can be
treated in a similar way.\hfill $\square $

\bigskip

\begin{lemma}
\label{EstOprel}Assume that $A1-$ $A2$ hold. For any $R\in \mathbb{Z}$, any $%
t\in (0,1)$ and $N\in \left\{ n-1,n\right\} $%
\begin{equation*}
\left\{ O_{\tau (nt)+R,N}|\tau (n)\geq nt\right\} \overset{d}{\rightarrow }%
\frac{1}{\Theta _{R}}-1.
\end{equation*}
\end{lemma}

\textbf{Proof}. We have 
\begin{equation*}
O_{\tau (nt)+R,N}=f_{\tau (nt)+R,N}(0)\beta _{N}(\tau (nt)+R)\left( \frac{%
e^{S_{\tau (nt)+R}}}{1-f_{0,\tau (nt)+R}(0)}-1\right) .
\end{equation*}%
Applying Lemmas \ref{Lpreelimit1}, \ref{Lprelim3}, and \ref{Lprelim4} we get%
\begin{equation*}
\left\{ O_{\tau (nt)+R,N}|\tau (n)\geq nt\right\} \overset{d}{\rightarrow }%
\frac{1}{\Theta _{R}}-1
\end{equation*}%
as desired.\hfill$\square$

\bigskip

\begin{lemma}
\label{DeltaVanish}Assume that $A1-$ $A2$ hold. Then%
\begin{equation*}
\Delta _{m,n}\overset{p}{\rightarrow }1\text{ as }\min \left( m,n-m\right)
\rightarrow \infty .
\end{equation*}
\end{lemma}

\textbf{Proof}. Clearly, it sufficies to show that%
\begin{equation}
\frac{1-f_{0,n}(0)}{1-f_{0,n-1}(0)}\overset{p}{\rightarrow }1\text{ as }%
n\rightarrow \infty .  \label{ratio1}
\end{equation}%
To verify this write%
\begin{equation*}
T_{\tau (n-1),n-1}^{\ast }:=\frac{1-f_{0,\tau (n-1)}(f_{\tau (n-1),n}(0))}{%
1-f_{0,\tau (n-1)}(f_{\tau (n-1),n-1}(0))}
\end{equation*}%
and%
\begin{equation*}
T_{l,n-l-1}:=\frac{1-f_{l,0}^{-}(f_{0,n-l-1}^{+}(f_{n-1}(0)))}{%
1-f_{l,0}^{-}(f_{0,n-l-1}^{+}(0))}.
\end{equation*}%
Since, evidently, there exists $s\in \lbrack 0,1)$ such that 
\begin{equation*}
\lim_{n-l\rightarrow \infty }f_{n-l-1}^{+}(s)\neq s,
\end{equation*}%
it follows from Lemma 1' and Theorem 5 in \cite{AK71} and (\ref{asq}) that $%
\mathbf{\hat{P}}^{+}$ a.s.%
\begin{equation*}
\lim_{n-l\rightarrow \infty
}f_{0,n-l-1}^{+}(f_{n-1}(0))=\lim_{n-l\rightarrow \infty
}f_{0,n-l-1}^{+}(0)=q^{+}.
\end{equation*}%
Hence we conclude that $\mathbf{\hat{P}}$ a.s. 
\begin{eqnarray*}
&&\lim_{\min (l.n-l)\rightarrow \infty }T_{l,n-l-1} \\
&=&\lim_{\min (l.n-l)\rightarrow \infty }\frac{%
1-f_{l,0}^{-}(f_{0,n-l-1}^{+}(f_{n-1}(0)))}{e^{S_{l}^{-}}}\times \frac{%
e^{S_{l}^{-}}}{1-f_{l,0}^{-}(f_{0,n-l-1}^{+}(0))}=1.
\end{eqnarray*}%
To finish the proof of (\ref{ratio1}) it remains to observe that 
\begin{equation*}
\lim_{n\rightarrow \infty }\mathbf{P}\left( \tau (n)\neq \tau (n-1)\right)
=0,
\end{equation*}%
to check that 
\begin{equation*}
\mathbf{E}\left[ T_{\tau (n-1),n-1}^{\ast }|\tau (n-1)=l\right] =\mathbf{E}%
\left[ T_{l,n-l-1}|A_{l,n-l-1}\right] ,
\end{equation*}%
and to apply Lemma \ref{Cfran}.\hfill $\square $

\bigskip

\textbf{Proof of Theorem} \ref{Tclose}. To prove the statement of the
theorem it suffices to combine Lemmas \ref{EstOprel}, \ref{EstO}, and \ref%
{DeltaVanish}.\hfill$\square$

\section{\protect\bigskip Proof of Theorem \protect\ref{Tfar}}

\begin{lemma}
\label{Lratio}Under Assumptions $A1-A2$,%
\begin{equation}
\frac{O_{nt,n-1}}{O_{nt,n}}\overset{p}{\rightarrow }1\text{ \ as \ }%
n\rightarrow \infty .  \label{Oratio}
\end{equation}
\end{lemma}

\textbf{Proof}. We have%
\begin{equation*}
\frac{O_{nt,n-1}}{O_{nt,n}}=\frac{f_{nt,n-1}(0)}{f_{nt,n}(0)}\times \frac{%
1-f_{0,n-1}(0)}{1-f_{0,n}(0)}\times \frac{1-f_{nt,n}(0)}{1-f_{nt,n-1}(0)}.
\end{equation*}%
Since the process is critical,%
\begin{equation}
f_{nt,n}(0)\overset{d}{=}f_{0,n(1-t)}(0)\rightarrow 1\text{ \ a.s.}
\label{c1}
\end{equation}%
as $n\rightarrow \infty .$ Similarly to (\ref{ratio1}) 
\begin{equation}
\frac{1-f_{nt,n}(0)}{1-f_{nt,n-1}(0)}\overset{d}{=}\frac{1-f_{0,n(1-t)}(0)}{%
1-f_{0,n(1-t)-1}(0)}\overset{d}{\rightarrow }1  \label{c3}
\end{equation}%
as $n\rightarrow \infty .$ Combining (\ref{ratio1}), (\ref{c1}), and (\ref%
{c3}) gives (\ref{Oratio}).\hfill$\square$

\bigskip

\begin{lemma}
\label{LOleftRight}Under Assumptions $A1$ $-$ $A2$, for any $t\in (0,1)$ as $%
n\rightarrow \infty $%
\begin{equation}
\left\{ O_{nt,n}\left( 1-f_{nt,n}(0)\right) \left\vert \tau (n)<nt\right.
\right\} \overset{d}{\rightarrow }1,  \label{di1}
\end{equation}%
and%
\begin{equation}
\left\{ O_{nt,n}\left( 1-f_{0,nt}(0)\right) e^{-S_{nt}}\left\vert \tau
(n)\geq nt\right. \right\} \overset{d}{\rightarrow }1.  \label{di2}
\end{equation}
\end{lemma}

\textbf{Proof}. To prove (\ref{di1})\ it suffices to observe that 
\begin{equation*}
O_{nt,n}\left( 1-f_{nt,n}(0)\right) =\frac{1}{\alpha _{n}(nt)}\times
f_{nt,n}(0)-f_{nt,n}(0)\beta _{n}(nt)\left( 1-f_{nt,n}(0)\right)
\end{equation*}%
and to apply Lemma \ref{min2} and (\ref{min42}), while to demonstrate (\ref%
{di2}) one should write 
\begin{equation*}
O_{nt,n}(1-f_{0,nt}(0))e^{-S_{nt}}=f_{nt,n}(0)\beta _{n}(nt)\left[ 1-\left(
1-f_{nt,n}(0))\right) e^{-S_{nt}}\right]
\end{equation*}%
and to use Lemma \ref{min4} and (\ref{min41}).\hfill$\square$

\bigskip

\begin{lemma}
\label{Comm1}Under Assumptions $A1-$\textbf{\ }$A2$, for any $t\in (0,1)$ as 
$n\rightarrow \infty $%
\begin{equation*}
\left\{ \frac{1-f_{0,nt}(0)}{e^{S_{\tau (nt)}}}\left\vert \tau (n)\geq
nt\right. \right\} \rightarrow \zeta
\end{equation*}%
and%
\begin{equation*}
\left\{ \frac{1-f_{nt,n}(0)}{e^{S_{\tau (nt,n)}-S_{nt}}}\left\vert \tau
(n)<nt\right. \right\} \rightarrow \zeta
\end{equation*}%
where $\zeta $ is defined in (\ref{aszet}).
\end{lemma}

\textbf{Proof}. Our arguments follow the same line of reasoning as the proof
of Lemma \ref{Lprelim4}. We have%
\begin{eqnarray*}
\lefteqn{\mathbf{E}\left[ \exp \left\{ -\lambda \frac{1-f_{0,nt}(0)}{%
e^{S_{\tau (nt)}}}\right\} ;\tau (n)\geq nt\right] } \\
&=&\mathbf{E}\left[ \exp \left\{ -\lambda \frac{1-f_{0,nt}(0)}{e^{S_{\tau
(nt)}}}\right\} \right]  \\
&&\hspace{2cm}-\mathbf{E}\left[ \exp \left\{ -\lambda \frac{1-f_{0,nt}(0)}{%
e^{S_{\tau (nt)}}}\right\} ;\tau (n)<nt\right] .
\end{eqnarray*}%
Introducing once again two independent environmental sequences $\left\{
f_{n}^{-},n\geq 0\right\} $ and $\left\{ f_{n}^{+},n\geq 0\right\} ,$ setting%
\begin{equation*}
T_{\tau (nt),nt}^{\ast }:=\exp \left\{ -\lambda \frac{1-f_{0,nt}(0)}{%
e^{S_{\tau (nt)}}}\right\} ,T_{l,nt-l}:=\exp \left\{ -\lambda \frac{%
1-f_{l,0}^{-}(f_{0,nt-l}^{+}(0))}{e^{S_{l}^{-}}}\right\} 
\end{equation*}%
and observing that by (\ref{aszet})%
\begin{equation*}
\frac{1-f_{l,0}^{-}(f_{0,nt-l}^{+}(0))}{e^{S_{l}^{-}}}\rightarrow \zeta
^{-}(q^{+})=\zeta 
\end{equation*}%
$\mathbf{\hat{P}}$ a.s. as $\min (l,nt-l)\rightarrow \infty ,$ \ we conclude
by Lemma \ref{Cfran} that for any $\lambda \in \left( 0,\infty \right) $%
\begin{equation*}
\mathbf{E}\left[ \exp \left\{ -\lambda \frac{1-f_{0,nt}(0)}{e^{S_{\tau (nt)}}%
}\right\} \right] \rightarrow \mathbf{\hat{E}}\left[ \exp \left\{ -\lambda
\zeta \right\} \right] 
\end{equation*}%
as $n\rightarrow \infty $. Further, 
\begin{eqnarray*}
\lefteqn{\mathbf{E}\left[ \exp \left\{ -\lambda \frac{1-f_{0,nt}(0)}{%
e^{S_{\tau (nt)}}}\right\} ;\tau (n)<nt\right] } \\
&=&\mathbf{E}\left[ \exp \left\{ -\lambda \frac{1-f_{0,nt}(0)}{e^{S_{\tau
(n)}}}\right\} ;\tau (n)<nt\right]  \\
&=&\sum_{j=0}^{nt-1}\mathbf{E}\left[ \exp \left\{ -\lambda \frac{%
1-f_{0,j}(f_{j,nt}(0))}{e^{S_{j}}}\right\} ;\tau (n)=j\right]  \\
&=&\sum_{j=0}^{nt-1}\mathbf{E}\left[ \exp \left\{ -\lambda \frac{%
1-f_{0,j}(f_{j,nt}(0))}{e^{S_{j}}}\right\} ;\min_{0\leq
r<j}S_{r}>S_{j},\min_{j+1\leq r\leq n}S_{r}\geq S_{j}\right]  \\
&=&\sum_{j=0}^{nt-1}\mathbf{E}\left[ \exp \left\{ -\lambda \frac{%
1-f_{j,0}^{-}(f_{0,nt-j}^{+}(0))}{e^{S_{j}^{-}}}\right\} ;A_{j,n-j}\right] 
\\
&=&\sum_{j=0}^{nt-1}\mathbf{E}\left[ \exp \left\{ -\lambda \frac{%
1-f_{j,0}^{-}(f_{0,nt-j}^{+}(0))}{e^{S_{j}^{-}}}\right\} \left\vert
A_{j,n-j}\right. \right] \mathbf{P}\left( \tau (n)=j\right) .
\end{eqnarray*}%
By\ (\ref{aszet}) 
\begin{equation*}
\frac{1-f_{j,0}^{-}(f_{0,nt-j}^{+}(0))}{e^{S_{j}^{-}}}\rightarrow \zeta 
\end{equation*}%
$\mathbf{\hat{P}}$ a.s. as $\min (j,nt-j)\rightarrow \infty $.  Hence we see
by Lemma \ref{Lfran} that 
\begin{equation*}
\lim_{\min (j,nt-j)\rightarrow \infty }\mathbf{E}\left[ \exp \left\{
-\lambda \frac{1-f_{j,0}^{-}(f_{0,nt-j}^{+}(0))}{e^{S_{j}^{-}}}\right\}
\left\vert A_{j,n-j}\right. \right] =\mathbf{\hat{E}}\left[ \exp \left\{
-\lambda \zeta \right\} \right] .
\end{equation*}%
This and the generalized arcsine law give%
\begin{eqnarray*}
\lefteqn{\lim_{n\rightarrow \infty }\mathbf{E}\left[ \exp \left\{ -\lambda 
\frac{1-f_{j,0}^{-}(f_{0,nt-j}^{+}(0))}{e^{S_{j}^{-}}}\right\} ;\tau (n)<nt%
\right] } \\
&=&\mathbf{\hat{E}}\left[ \exp \left\{ -\lambda \zeta \right\} \right] 
\mathbf{P}\left( \tau \leq t\right) .
\end{eqnarray*}%
Thus,%
\begin{eqnarray*}
\lefteqn{\lim_{n\rightarrow \infty }\mathbf{E}\left[ \exp \left\{ -\lambda 
\frac{1-f_{0,nt}(0)}{e^{S_{\tau (nt)}}}\right\} ;\tau (n)\geq nt\right] } \\
&=&\mathbf{\hat{E}}\left[ \exp \left\{ -\lambda \zeta \right\} \right] 
\mathbf{P}\left( \tau \geq t\right) 
\end{eqnarray*}%
and the first statement of the lemma follows. The second statement can be
checked in a similar way. \hfill $\square $

\begin{lemma}
\label{LOinf} Under Assumptions $A1-$ $A2$, for any $t\in (0,1)$%
\begin{equation*}
O_{nt,n}\overset{d}{\rightarrow }\infty \text{ \ as }n\rightarrow \infty .
\end{equation*}
\end{lemma}

\textbf{Proof}. This statement is a direct corollary of Lemmas \ref%
{LOleftRight} and \ref{Comm1}.\hfill \hfill $\square $

\bigskip

\textbf{Proof of Theorem \ref{Tfar}. }Using Lemma \ref{BasicINeq} we have 
\begin{eqnarray*}
\frac{1}{1+O_{nt,n-1}(1-e^{-\lambda /O_{nt,n-1}})}\Delta _{nt,n}^{-1} &\leq &%
\mathbf{E}_{\mathbf{k}}\left[ \exp \left\{ -\lambda \frac{Z_{nt}}{O_{nt,n}}%
\right\} \left\vert T=n\right. \right] \\
&\leq &\frac{1}{1+O_{nt,n}(1-e^{-\lambda /O_{nt,n}})}\Delta _{nt,n}.
\end{eqnarray*}%
Now to complete the proof of the theorem it remains to observe that 
\begin{equation*}
\lim_{n\rightarrow \infty }\Delta _{nt,n}\overset{p}{=}1
\end{equation*}%
by Lemmas \ref{Lratio} and \ref{LOinf} and (\ref{c3}), and that by Lemma \ref%
{LOinf} 
\begin{equation*}
\lim_{n\rightarrow \infty }O_{nt,n}(1-e^{-\lambda /O_{nt,n}})\overset{d}{=}%
\lim_{n\rightarrow \infty }O_{nt,n-1}(1-e^{-\lambda /O_{nt,n}})\overset{d}{=}%
\lambda .
\end{equation*}

\begin{remark}
\label{roughlyspeaking} It follows from Lemmas \ref{LOleftRight} and \ref%
{Comm1} that, roughly speaking, 
\begin{equation*}
O_{nt,n}I\left\{ \tau (n)\geq nt\right\} \approx \frac{e^{S_{nt}}}{%
1-f_{0,nt}(0)}I\left\{ \tau (n)\geq nt\right\} \asymp e^{S_{nt}-S_{\tau
(nt)}}I\left\{ \tau (n)\geq nt\right\} 
\end{equation*}%
and%
\begin{equation*}
O_{nt,n}I\left\{ \tau (n)<nt\right\} \approx \frac{1}{1-f_{nt,n}(0)}I\left\{
\tau (n)<nt\right\} \asymp e^{S_{nt}-S_{\tau (nt,n)}}I\left\{ \tau
(n)<nt\right\} .
\end{equation*}%
Thus, given $\left\{ T=n\right\} $ the growth rate of $Z_{nt}$ depends
essentially on the location of the point $\tau (n)$ of the global maximum of 
$\left\{ S_{k},0\leq k\leq n\right\} $  with respect to the moment $nt.$ If $%
\tau (n)\geq nt$ then this growth rate is specified by the past local
minimum of the associated random walk while if $\tau (n)<nt$ then it depends
on the value of this random walk at the point of prospective minimum of its
remaining piece.
\end{remark}

\section*{Acknowledgment}

Part of this work was carried out whilst V.A.V. was visiting A.E.K. at The
Department of Mathematical Sciences at Heriot Watt University during the
months of June and July 2006. V.A.V. is grateful for the department's
hospitality during this period.

\end{document}